\documentclass{ifacconf}

\usepackage{graphicx}      % include this line if your document contains figures
\usepackage{natbib}        % required for bibliography
\usepackage[dvips]{epsfig}    % or this line, depending on which
\usepackage{xcolor}

\usepackage{amssymb}
\usepackage{amsmath}          % Include this line if your
\usepackage{enumerate}

{\theoremheaderfont{\upshape} \newtheorem{myprop}{Proposition}}

\begin{document}
\begin{frontmatter}

\title{Modified global finite-time quasi-continuous second-order robust feedback control}

%\thanks[]{\textcolor[rgb]{0.00,0.00,1.00}{Authors accepted manuscript}}

\author[1]{Michael Ruderman},
\author[2]{Denis Efimov}

\address[1]{University of Agder, Department of Engineering Sciences, \\ Grimstad, Norway  \\
Email (correspondence): \tt\small michael.ruderman@uia.no}
\address[2]{INRIA, Univ. Lille, CNRS, CRIStAL, Lille, France.}

%%%%%%%%%%%%%%%%%%%%%%%%%%%%%%%%%%%%%%%%%%%%%%%%%%%%%%%%%%%%%%%%%%%%%%%%%%%%%%%%
\begin{abstract}
A non-overshooting quasi-continuous sliding mode control with sub-optimal damping was recently introduced in \cite{ruderman2024robust} for perturbed second-order systems. The present work proposes an essential modification of the nonlinear control law which (i) allows for a parameterizable control amplitude limitation in a large subset of the initial values, (ii) admits an entire state-space $\mathbb{R}^2$ (that was not given in \cite{ruderman2024robust}) for the finite-time control, and finally (iii) enables for the found analytic solution of the state trajectories in the unperturbed case. The latter allows also for an exact estimation of the finite convergence time, and open an avenue for other potentially interesting analysis of the control properties in the future. For a perturbed case, the solution-based and Lyapunov function-based approaches are developed to show the uniform global asymptotic stability. The proposed robustness and convergence analysis are accompanied by several illustrative numerical examples.  
\end{abstract}

\begin{keyword}
Second-order dynamics \sep robust control \sep finite-time convergence \sep quasi-continuous control \sep non-overshooting \sep sub-optimal damping \sep analytic solutions
\end{keyword}

\end{frontmatter}

%%%%%%%%%%%%%%%%%%%%%%%%%%%%%%%%%%%%%%%%%%%%%%%%%%%%%%%%%%%%%%%%%%%%%%%%%%%%%%%%
\section{Introduction}
\label{sec:1}

Design of robust control laws for second order systems in the canonical forms is an important area of research, see e.g. \cite{Tsypkin1984relay,Kwakernaak1993,aastrom2006}, since a wide spectrum of mechanical and electrical systems are described by this class of models. There are different canonical representations of the second order systems, and among them the double integrator dynamics are very popular.

The diversity of the existing control solutions for a disturbed double integrator model is related with different quality requirements that appear in applications. Among them, it is worth to mention the convergence rate of the regulation error, tolerance to the perturbations, the complexity of the control algorithm and its simplicity for the implementation (\cite{khalil2002}). A group of methods that answer all these requirements is developed in the sliding mode control theory (\cite{shtessel2014,utkin2020}). These controls frequently provide the finite-time convergence to zero of the trajectories in the closed-loop system, uniformly in the matched and properly bounded disturbances, with simple structure of the control law.

An interesting higher order sliding mode control algorithm (see \cite{Levant2005} for the introduction of this kind of stabilizers) was proposed recently in \cite{ruderman2024robust}, which in addition to other positive characteristics mentioned above also guarantees the absence of overshooting in the presence of disturbances. This property is important for safety critical systems (\cite{Polyakov2023}), where it is required to ensure that the system does not visit certain domains of the state space during the transients. However, the control from \cite{ruderman2024robust} also possesses several drawbacks: first, it was not defined for the whole state space of the system, second, the control amplitude was quickly growing close to the origin in the first and third quadrants of the state space. In order to solve these issues, a new modified control law is proposed in this paper, which is a uniformly bounded in the first and third quadrants of the state space, well defined for all states in $\mathbb{R}^2$, and moreover, the solutions of the closed-loop nonlinear system can be analytically calculated for a dominant part of the state space. The latter achievement allows for the convergence time to be accurately evaluated and for further investigations. The finite-time stability of the closed-loop system is proven using the Lyapunov function method. Several comparative simulations illustrate also the properties and advantages of the proposed control algorithm.

%%%%%%%%%%%%%%%%%%%%%%%%%%%%%%%%%%%%%%%%%%%%%%%%%%%%%%%%%%%%%%%%%%%%%%%%%%%%%%%%
\section{Control system}
\label{sec:2}

We consider a generic class of perturbed second-order dynamic systems described by 
\begin{eqnarray}
\nonumber \dot{x}_1(t) &=& x_2(t), \\
\dot{x}_2(t) &=& u(t) + d(t),
\label{eq:2:1}
\end{eqnarray}
where the unknown disturbance value $d(t)\in\mathbb{R}$ is matched with the
control quantity $u(t)$, and is (Lebesgue) measurable,
while $\|d\|_{\infty}=\text{ess}\sup_{t\geq0}|d(t)|\leq D$. The dynamic  
state variables are assumed to be available for a feedback control and 
denoted in the vector form (for the sake of a brief notation in course of analysis) by $x(t)=\left(x_{1}(t),x_{2}(t)\right)^{\top} \in \mathbb{R}^{2}$. The upper-bound $0 < D < \infty$ of dynamic disturbances is assumed to be known. Further we note that the inertial and input-gaining terms, in case of a physical plant model equivalent to \eqref{eq:2:1}, can be correspondingly incorporated into the control parameters when designing $u(t)$.  

\begin{figure}[!h] \centering
\includegraphics[width=0.98\columnwidth]{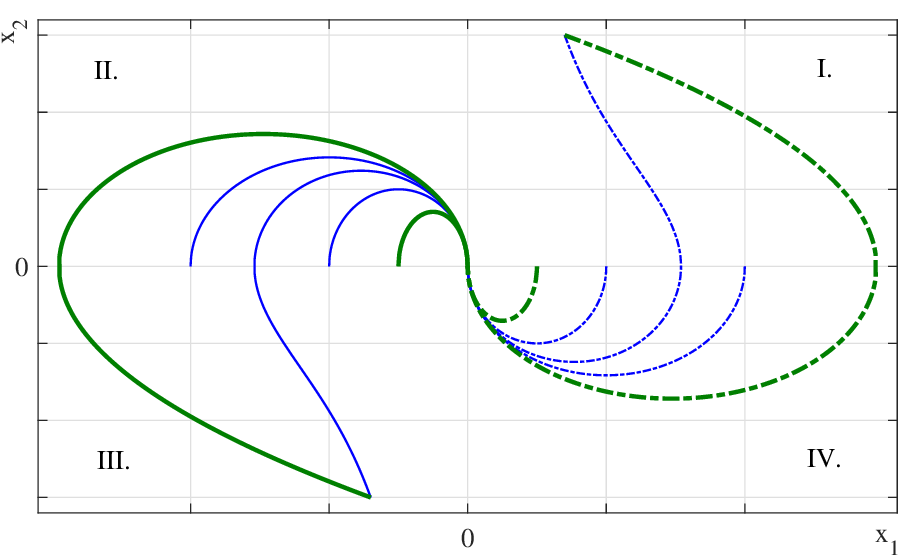}
\caption{Phase portrait of the unperturbed (i.e. $d=0$) second-order control system: control (\cite{ruderman2024robust}) drawn by the thin-blue lines, and control \eqref{eq:2:2},\eqref{eq:2:3} drawn by the thick-green lines.}
\label{fig:2:1}
\end{figure}

In \cite{ruderman2024robust}, a novel quasi-continuous robust control $u(t)$ with a fast nonlinear damping law (introduced in \cite{ruderman2021a} and detailed with experiments in \cite{ruderman2022}) was proposed for \eqref{eq:2:1}. The existence of solutions, uniform asymptotic stability, finite time convergence, and robustness against the disturbances
$\|d\|_{\infty}\leq D$ of the closed-loop system were shown with an accompanying
analysis in \cite{ruderman2024robust}. Note that this control law was only defined on the invariant for the closed-loop system set $\{x\in\mathbb{R}^{2}:x_1\ne0\}\cup\{0\}$. The phase portrait of that control system is qualitatively  shown in Fig. \ref{fig:2:1} by the thin blue lines. Note that the depicted notation of the quadrants I.-IV. in the phase-plane will further be used in analysis of the control system.

The proposed modification of the quasi-continuous robust control (\cite{ruderman2024robust}) is given by
\begin{eqnarray}
\label{eq:2:2}
u(t)  & = & -\gamma \, \textrm{sign}\bigl(x_1(t)\bigr) - \delta(t), \\[1mm]
\delta(t) &=& \left\{
             \begin{array}{ll}
             \bigl|x_2(t)\bigr| \bigl|x_1(t)\bigr|^{-1} x_2(t) & \hbox{ if } \; x_1(t) x_2(t) < 0,  \\[1mm]
              0 & \hbox{ else},
             \end{array}
           \right.
\label{eq:2:3}
\end{eqnarray}
where $\gamma > D$ is the (same) single design parameter. Also the nonlinear damping term when $\delta \neq 0$ has the same analytic form as in \cite[eq.~(2)]{ruderman2024robust}. The main principal difference is that no damping is effective in the I. and III. quadrants, cf. \eqref{eq:2:3}. This way, a fastest possible transition into the IV., respectively II., quadrant is achieved if the control amplitude is limited by $\gamma$ when $x(0) \in$ I. $\vee$ III. quadrat. The same $\max|u| = \gamma$ can also be guaranteed if $x(0) \in$ II. $\vee$ IV. quadrat and, at the same time, the initial conditions are sufficiently small (e.g., $x(0) \in C_a$, cf. \eqref{eq:3:2:1} below). Recall that the convergence of the state trajectories to the globally stable zero equilibrium occurs always within either II. or IV. quadrant, cf. Fig. \ref{fig:2:1}. The convergence properties and the analytic solution of the control system \eqref{eq:2:1}--\eqref{eq:2:3} will be provided below in sections \ref{sec:3} and \ref{sec:4}, while at once the following remarks are worthily making.

\begin{rem}
\label{rem:1} Directly recognizable (by comparing the control equations \eqref{eq:2:2},\eqref{eq:2:3} and \cite[eq.~(2)]{ruderman2024robust} is that the state trajectories differ only in the I. and III. quadrants. Otherwise, in the II. and IV. quadrants, the state trajectories of the control systems (1)--(3) and (\cite{ruderman2024robust}) coincide with each other for the same initial conditions $x(0) \in \{x\in\mathbb{R}^{2} \,: \; x_1 x_2 < 0 \}$.
\end{rem}

\begin{rem}
\label{rem:2} Unlike for the control system in \cite{ruderman2024robust}, where 
the admissible set of the initial conditions $x(0)\in\{x\in\mathbb{R}^{2}\,:\,x_{1}\neq0\}\cup\{0\}$ excludes the
$x_2$-axis if $x_2 \neq 0$, which is reasonable for multiple but not all possible application scenarios, the initial conditions for the control system \eqref{eq:2:1}--\eqref{eq:2:3} satisfy $x(0) \in \mathbb{R}^{2}$. %\mathbb{X}=
\end{rem}

\begin{rem}
\label{rem:3} 
The length of the trajectories within the I. and III. quadrants, and so does the time of reaching $x_2 = 0$, increases comparing to the control solution (\cite{ruderman2024robust}). However, for any $x(0) \in \{x\in\mathbb{R}^{2} \,: \; x_1 x_2 \geq 0 \}$, the solutions of \eqref{eq:2:1}--\eqref{eq:2:3} are the fastest possible within the I. and III. quadrants, in case the control amplitude has to be limited, i.e., $\bigl|u(t)\bigr| \leq \gamma$.
\end{rem}

%%%%%%%%%%%%%%%%%%%%%%%%%%%%%%%%%%%%%%%%%%%%%%%%%%%%%%%%%%%%%%%%%%%%%%%%%%%%%%%%
\section{Unperturbed behavior}
\label{sec:3}

\begin{figure}[!h] \centering
\includegraphics[width=0.98\columnwidth]{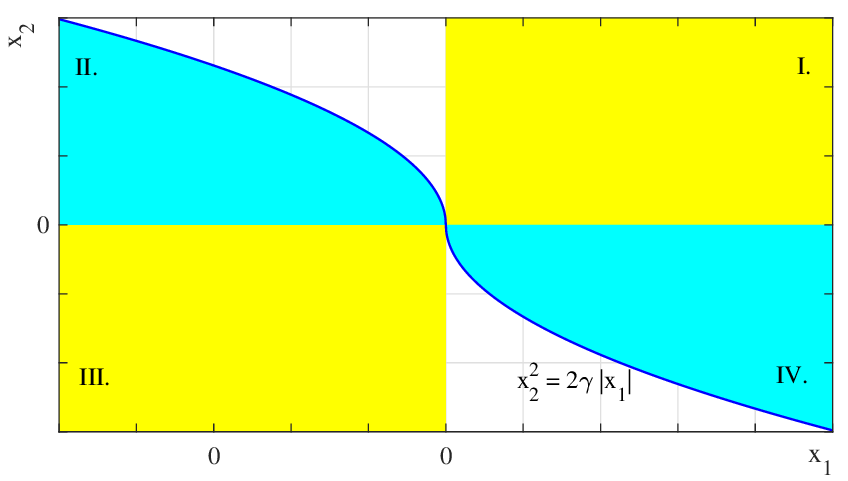}
\caption{State-space of the control system \eqref{eq:2:1}--\eqref{eq:2:3} with the regions of 
unperturbed (i.e. $d=0$) analytic solutions for $U$ (yellow shadowed) and for $C_a$ (cyan shadowed).}
\label{fig:3:1}
\end{figure}

In the following, we analyze first the behavior of an unperturbed control system \eqref{eq:2:1}--\eqref{eq:2:3}, i.e. for $d=0$, and that separately for I. and III. and then for II. and IV. quadrants. More rigorously, we define the quadrants as 
\begin{eqnarray}
\label{eq:3:1}
  x \; &  \in & \; \textrm{I.} \cup \textrm{III.} \equiv U = \{x\in\mathbb{R}^{2} \,: \; x_1 x_2 \geq 0 \}, \\
  x \; &  \in &  \; \textrm{II.} \cup \textrm{IV.} \equiv C = \{x\in\mathbb{R}^{2} \,: \; x_1 x_2 < 0 \},
\label{eq:3:2}
\end{eqnarray}
respectively, where $U$ refers to an undamped subset and $C$ to a convergence subset of trajectories in the state-space. Recall that due to $\dot{x}_1 = x_2$ and, thus, clockwise evolvement of the state trajectories, the convergence towards zero equilibrium appears only when $x(t) \in C$, cf. Fig. \ref{fig:2:1}. Also it is worth noting that zero equilibrium itself belongs to $U$, so that a control action with zero damping, cf. \eqref{eq:2:3}, appears also at $x(t) = (0,0)$, namely in the sliding-mode. Respectively, no control singularities can occur at $x_2$-axis if $x_2 \neq 0$, cf. \cite{ruderman2024robust}. For $U$, we provide the analytic solutions in section \ref{sec:3:sub:1}. For $C$, the analytic solution (section \ref{sec:3:sub:2}) is obtained for the subset 
\begin{equation}\label{eq:3:2:1}
C_a = \bigl\{ x\in C \; : \;  x_2^2 < 2 \gamma |x_1| \bigr\}
\end{equation}
of initial conditions, see Fig. \ref{fig:3:1} (cyan shadowed regions). 

Due to an obvious symmetry of the trajectory solutions and, therefore, without loss of generality, the solutions will be developed below in the I. quadrant for $U$ and in the IV. quadrant for $C_a$, cf. the dashed lines in Fig. \ref{fig:2:1}.

Further we note that since a transition from $U$ to $C_a$ is possible at $x_2 = 0$ only, no control jumps appear when switching $\delta$ in the control law \eqref{eq:2:2}, \eqref{eq:2:3}. Respectively, the state trajectories and so the control value remain smooth outside zero equilibrium.

\subsection{Solution in I. (respectively III.) quadrant}
\label{sec:3:sub:1}

For an arbitrary initial state\footnote{Note that since the control law \eqref{eq:2:2}, \eqref{eq:2:3} yields the closed-loop system with a switching dynamics, any time instant where a state trajectory reaches the $x_1$- or/and $x_2$-axis and, therefore, attains $U$ can be seen as an initial time $t=0$.} $x(0)=\bigl(x_{1}(0),x_{2}(0)\bigr)^{\top} \geq 0$, and the control system \eqref{eq:2:1}--\eqref{eq:2:3} with $d=0$, the analytic solution of a double-integrator with the non-homogenous part $u(t) = -\gamma$, cf. \eqref{eq:2:3}, is given by, cf. \cite{Athans63},
\begin{eqnarray}
\label{eq:3:3}
  x_1(t) &=& x_1(0) + x_2(0) \, t - \dfrac{\gamma \, t^2}{2}, \\[1mm]
  x_2(t) &=& x_2(0) - \gamma \, t 
\label{eq:3:4}
\end{eqnarray}
for the instants of time $t\geq0$ when the trajectory stays in $U$. One can recognize that the analytic solution \eqref{eq:3:3}, \eqref{eq:3:4} results in a parabolic trajectory, cf. Fig. \ref{fig:2:1}, which starts at $\bigl(x_1(0), x_2(0)\bigr)$ and ends when $x_2(T_{U_{+}}) = 0$ for the finite reaching time 
\begin{equation*}
T_{U_{+}} = \dfrac{x_2(0)}{\gamma}.
\end{equation*}
The latter is obtained for the I. quadrant by solving the final-value problem of \eqref{eq:3:4} with $x_2(t)=0$. Easy to recognize is that the transient time for the entire $U$ is, respectively,  
\begin{equation}\label{eq:3:6}
T_{U_{\pm}} = \dfrac{\bigl|x_2(0)\bigr|}{\gamma}.
\end{equation}

\vspace{2mm}

\begin{rem}
\label{rem:4}
It is worth noting that also zero solution of the unperturbed system \eqref{eq:2:1}--\eqref{eq:2:3} for $x(0) = (0,0)$ can equally be proven by using \eqref{eq:3:3}, \eqref{eq:3:4}. Indeed, depending on whether in the I. or III. quadrant, the solution \eqref{eq:3:4} with $x_2(0)=0$ transforms to the set of equations
\begin{equation}\label{eq:3:7}
x_2(t) = -\gamma \, t \quad \hbox{ and } \quad x_2(t) = \gamma \, t,
\end{equation}
respectively. If $x_1(0)\ne0$ this leads to an instantaneous exit of the trajectory to the domain $C$. Otherwise, for $x_1(0)=0$, the system is reduced to sliding mode dynamics (it corresponds to $\delta = 0$ of the control law \eqref{eq:2:2}, \eqref{eq:2:3}) $$\dot{x}_1(t)=x_2(t),\; \dot{x}_2(t)=-\gamma\,\text{sign}(x_1(t))$$ that has the equilibrium at the origin, and for all $t > 0$ when switching happens between $\pm \gamma$ with the infinite frequency so that $x(t)$ remains zero. Consequently, $x(t) = (0,0)$ is the unique solution of \eqref{eq:2:1}--\eqref{eq:2:3} for all $t > 0$ if $d=0$ and $x(0) = (0,0)$.
\end{rem}

\subsection{Solution in IV. (respectively II.) quadrant $x(0) \in C_a$}
\label{sec:3:sub:2}

In the following, we assume $x(0) \in C_a$ while considering $x(0)$ to be in the IV. quadrant (see Fig. \ref{fig:3:1}). This is without loss of generality since the solutions in the IV. and II. 
quadrants are symmetric to each other about the origin. The solution for $x(0) \in C_a$ in the II. quadrant can then be derived in the similar manner.

\begin{myprop}
\label{prop:1} The homogeneous solution of the feedback controlled system
\eqref{eq:2:1}--\eqref{eq:2:3} with $d=0$, $x(0) \in C_a$ and $x_1(0) > 0$ is given by
\begin{eqnarray}
\label{eq:3:1:1}
  x_1(t) &=& - \frac{\gamma}{\omega^2} \cos (\omega t + \phi) + B, \\[1mm]
  x_2(t) &=& \frac{\gamma}{\omega} \sin (\omega t + \phi),
\label{eq:3:1:2}
\end{eqnarray}
where
\begin{equation}\label{eq:3:1:3}
B = \dfrac{\gamma \, x_1^2(0)}{2 \gamma \, x_1(0) - x_2^2(0)}, % this was typo in the numerator, SW implementation was using correctly x_1^2(0) in numerator
\end{equation}
and the harmonic parameters are
% my harmonic parameters look OK, I'm going to check your math version and share an evaluating script
\begin{equation}
\label{eq:3:1:4}
  \phi = \pi + \arccos \biggl(1 - \dfrac{x_2^2(0)}{\gamma \, x_1(0)} \biggr), \quad   \omega = \sqrt{\dfrac{\gamma}{B}}. 
\end{equation}
\end{myprop}

\vspace{2mm}

\begin{pf}
First, consider the proposed solution \eqref{eq:3:1:1}, \eqref{eq:3:1:2} with the unknown values of $B$ and harmonic parameters. The property $\dot{x}_1(t)=x_2(t)$ can be verified by direct differentiation. Taking then the time derivative of \eqref{eq:3:1:2}, one obtains 
\begin{equation}\label{eq:3:1:5}
\dot{x}_2(t) = u(t) = \gamma \cos(\omega t + \phi).
\end{equation}
Given the finite-time ($0 < T_{C_a} < \infty$) convergence of the control \eqref{eq:2:2}, \eqref{eq:2:3} in the IV. quadrant, see \cite{ruderman2024robust} for details and proofs, the control value at $T_{C_a}$ (i.e. in the state-space origin) satisfies  
\begin{equation}\label{eq:3:1:6}
\dot{x}_2(T_{C_a}) = \gamma = \gamma \cos(\omega T_{C_a} + \phi).
\end{equation}
Using the right-hand-side equality in \eqref{eq:3:1:6} and solving for both equations \eqref{eq:3:1:1} and \eqref{eq:3:1:2} the initial value (with the given $x(0)$) and final value  $x(T_{C_a}) = 0$ problems, results in the parameter values \eqref{eq:3:1:3} and \eqref{eq:3:1:4}.

Next, it is evident that the unperturbed closed-loop system \eqref{eq:2:1}-\eqref{eq:2:3} in the IV. quadrant transforms into   
\begin{equation}\label{eq:3:1:7}
\dot{x}_2(t) - \dfrac{x_2^2(t)}{x_1(t)} + \gamma = 0.
\end{equation}
Substituting \eqref{eq:3:1:1}, \eqref{eq:3:1:2}, and \eqref{eq:3:1:5} with \eqref{eq:3:1:3} and \eqref{eq:3:1:4} into \eqref{eq:3:1:7}, yields that the differential equation \eqref{eq:3:1:7}, and so \eqref{eq:2:1}-\eqref{eq:2:3} for the IV. quadrant, holds. This completes the proof.
\end{pf}

\vspace{2mm}

\begin{rem}
\label{rem:5} The finite time of the convergence to origin is then obtained by solving directly \eqref{eq:3:1:6} with \eqref{eq:3:1:4} as 
\begin{equation}\label{eq:3:1:8}
T_{C_a} = x_1(0) \, \dfrac{\pi - \arccos \biggl(1 - \dfrac{x_2^2(0)}{\gamma x_1(0)} \biggr)}{\sqrt{\gamma x_1(0) - x_2^2(0)}}.
\end{equation}
Note that \eqref{eq:3:1:8} is computed for the first period of the otherwise periodic solution of \eqref{eq:3:1:6} with respect to $T_{C_a}$.
\end{rem}

\subsection{Numerical illustration}
\label{sec:3:sub:3}

\begin{figure}[!h] \centering
\includegraphics[width=0.98\columnwidth]{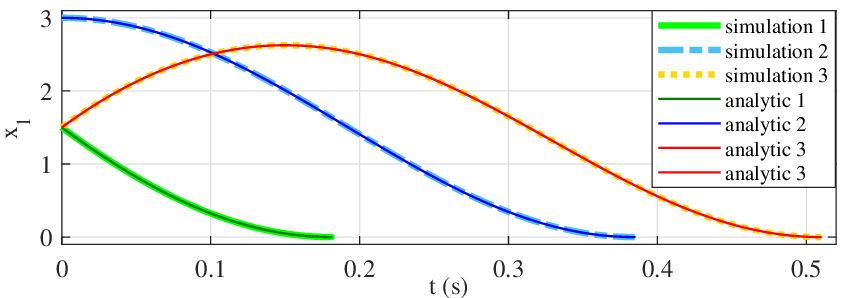}
\includegraphics[width=0.98\columnwidth]{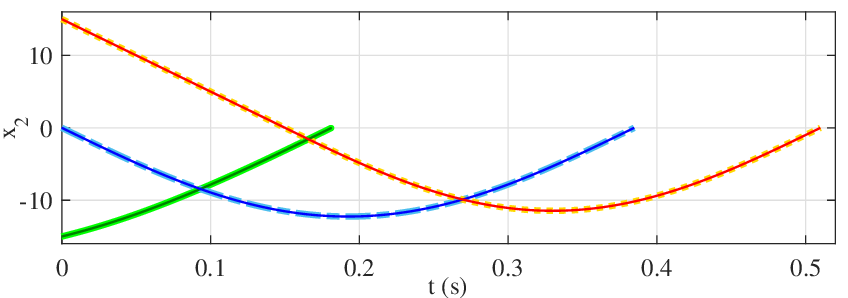}
\includegraphics[width=0.98\columnwidth]{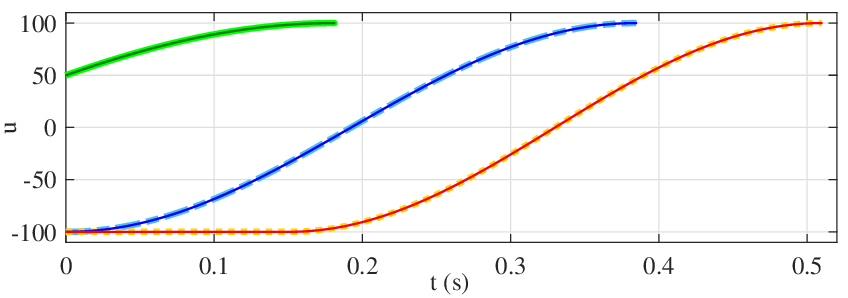}
\caption{Comparison of analytic solutions and numerical simulations of the unperturbed ($d=0$) closed-loop system \eqref{eq:2:1}, \eqref{eq:2:2} for three different initial conditions.}
\label{fig:3:2}
\end{figure}

The obtained analytic solutions \eqref{eq:3:3}, \eqref{eq:3:4} and \eqref{eq:3:1:1}--\eqref{eq:3:1:4} are compared with numerical simulation of the closed-loop system \eqref{eq:2:1}, \eqref{eq:2:2} for the case $d=0$ and three different $x(0) \in U \vee C_a$. The assign feedback gain is $\gamma = 100$. A standard fixed-step numerical solver of Matlab/Simulink is used with the assigned sampling rate of  10 kHz. Also the control values $u(t)$ for both, the analytic solutions and numerical simulation, are compared for the sake of completeness. The results are depicted in Fig. \ref{fig:3:2}.

%%%%%%%%%%%%%%%%%%%%%%%%%%%%%%%%%%%%%%%%%%%%%%%%%%%%%%%%%%%%%%%%%%%%%%%%%%%%%%%%

\section{Perturbed behavior}
\label{sec:4}

Note that for $\|d\|_{\infty}\leq D$ and $\gamma>D$, the system
(\ref{eq:2:1})--(\ref{eq:2:3}) admits the unique equilibrium at
the origin, which can be shown by analysis of the equation $\dot{x}=0$
(a well-known result in the sliding mode control theory obtained for
the twisting control).

For the disturbance-free case, the closed-loop system stability can
be shown using an energy-like Lyapunov function
\[
E(x)=\gamma|x_{1}|+\frac{1}{2}x_{2}^{2},\;\dot{E}(x)\leq0,
\]
which, however, cannot be used if $d\ne0$, since its derivative is
not sign-definite. For the quasi-continuous control of (\cite{ruderman2024robust}),
the uniform asymptotic stability of the closed loop system was shown
provided
\begin{equation}
\gamma>D^{1.5}+D+\frac{1}{2},\label{eq:gamma_old}
\end{equation}
and analyzing the following Lyapunov function candidate:
\begin{align}
V(x) & =\frac{1}{2}\xi^{\top}(x)\left[\begin{array}{cc}
2\gamma & \varepsilon\\
\varepsilon & 1
\end{array}\right]\xi(x),\label{eq:V_old}\\
\xi(x) & =\left[\begin{array}{c}
\sqrt{|x_{1}|}\text{sign}(x_{1})\\
x_{2}
\end{array}\right],\nonumber 
\end{align}
which is positive definite providing that the parameter $\varepsilon\in(0,\frac{2}{3})$,
and continuously differentiable for all $x\in\mathbb{R}^{2}\setminus\{x\in\mathbb{R}^{2}:x_{1}=0\}$.
Unfortunately, this $V(x)$ cannot be used for the analysis of (\ref{eq:2:1})--(\ref{eq:2:3})
since on the set $U$ the damping is suppressed, in comparison
to (\cite{ruderman2024robust}); then the system is just stable while
$x_{2}(t)x_{1}(t)\geq0$.

For demonstrating the uniform global asymptotic stability
of (\ref{eq:2:1})--(\ref{eq:2:3}), the authors found two approaches. 

\subsection{Solution-based approach}

Note that for $x_{2}(t)x_{1}(t)<0$, on the domain $C$, the system
(\ref{eq:2:1})--(\ref{eq:2:3}) coincides with the one studied in
(\cite{ruderman2024robust}), where the following properties have
been established for the respective closed-loop dynamics:
\begin{eqnarray*}
\dot{x}_{1}(t) & = & x_{2}(t),\\
\dot{x}_{2}(t) & = & -\gamma\,\textrm{sign}\bigl(x_{1}(t)\bigr)-\frac{\bigl|x_{2}(t)\bigr|}{\bigl|x_{1}(t)\bigr|}x_{2}(t)+d(t),
\end{eqnarray*}

\begin{itemize}
\item the solutions of this system are well-defined for all $t\geq0$ and
initial conditions $x(0)\in C\cup\{0\}$ with $\|d\|_{\infty}\leq D$;
\item the set $C\cup\{0\}$ is forward invariant for any $\|d\|_{\infty}\leq D$;
\item the origin is uniformly finite-time attractive and stable for the
initial conditions $x(0)\in C\cup\{0\}$ and the disturbances $\|d\|_{\infty}\leq D$.
\end{itemize}
The principal intuition behinds these properties is that under (\ref{eq:gamma_old})
the first term is always dominating the disturbance, then (\ref{eq:V_old})
can be applied for the analysis.

Therefore, it is left to show that for $x(0)\in U\setminus\{0\}$
and $\|d\|_{\infty}\leq D$, any respective trajectory of the system
leaves $U$ and goes to $C\cup\{0\}$ in a finite time, then uniform
global finite-time attractiveness of the origin is substantiated.
In the set $U$ the system (\ref{eq:2:1})--(\ref{eq:2:3}) takes
the form
\begin{eqnarray*}
\dot{x}_{1}(t) & = & x_{2}(t),\\
\dot{x}_{2}(t) & = & -\gamma\,\textrm{sign}\bigl(x_{1}(t)\bigr)+d(t).
\end{eqnarray*}
For $\gamma>D$, this system can stay on the axes $x_{1}=0$ or $x_{2}=0$
out the origin during an isolated instant of time only (in both cases
$\dot{x}_{1}$ and $\dot{x}_{2}$ are separated with zero, respectively,
that forces the system to enter in the interior of the set $U$).
Then, using the propriety that $x_{2}(t)x_{1}(t)>0$ we get
\begin{eqnarray*}
\frac{d|x_{1}(t)|}{dt} & = & |x_{2}(t)|,\\
\frac{d|x_{2}(t)|}{dt} & = & -\gamma+\text{sign}(x_{2}(t))d(t),
\end{eqnarray*}
and the following estimates in the time domain follow:
\begin{gather*}
|x_{2}(0)|-(\gamma+D)t\\
\leq|x_{2}(t)|=|x_{2}(0)|-\gamma t+\int_{0}^{t}\text{sign}(x_{2}(s))d(s)ds\\
\leq|x_{2}(0)|-(\gamma-D)t
\end{gather*}
which is satisfied for $t\in[0,T_{U}]$, where the finite time of
leaving the set $U$ has the bounds:
\[
\frac{|x_{2}(0)|}{\gamma+D}\leq T_{U}\leq\frac{|x_{2}(0)|}{\gamma-D}.
\]
Consequently, uniform (in the sense of being independent in $d$ with
$\|d\|_{\infty}\leq D$) global finite-time attractiveness of the
origin is proven for (\ref{eq:2:1})--(\ref{eq:2:3}).

In order to show the stability of the origin, we need anyway to find
a suitable Lyapunov function, which brings us to another approach
for substantiation of convergence and stability of this system.

\subsection{Lyapunov function approach}

The main result of this section is formulated below:
\begin{thm}
The origin for (\ref{eq:2:1})--(\ref{eq:2:3}) is uniformly globally
asymptotically stable if 
\begin{equation}
\gamma>2\sqrt{2}D^{1.5}+D+\frac{1}{2}.\label{eq:gamma_new}
\end{equation}
\end{thm}
\begin{pf}
Let us consider the following Lyapunov function candidate for $\gamma>D+\eta$
with some parameter $\eta>0$:
\begin{align}
V(x) & =\frac{1}{2}\xi^{\top}(x)\left[\begin{array}{cc}
2(\gamma-D-\eta) & \epsilon(x)\\
\epsilon(x) & 1
\end{array}\right]\xi(x),\label{eq:V_new}\\
\xi(x) & =\left[\begin{array}{c}
\sqrt{|x_{1}|}\text{sign}(x_{1})\\
x_{2}
\end{array}\right],\;\epsilon(x)=\varepsilon\begin{cases}
1 & \text{if }x_{1}x_{2}<0,\\
0 & \text{otherwise,}
\end{cases}\nonumber
\end{align}
which is positive definite providing that the parameter $\varepsilon\in(0,\sqrt{2(\gamma-D-\eta)})$.
The function $V$ is continuous on $\mathbb{R}^{2}$ since the coupling
term $\sqrt{|x_{1}|}\text{sign}(x_{1})x_{2}$ vanishes at the set
$\mathcal{X}$, where $\mathcal{X}=\{x\in\mathbb{R}^{2}:x_{1}=0\vee x_{2}=0\}$,
and it is continuously differentiable for all $x\in\mathbb{R}^{2}\setminus\mathcal{X}$.
As it has been mentioned in the previous subsection, the system can
stay in the set $\mathcal{X}\setminus\{0\}$ at an isolated time instant
only, and the origin is the equilibrium. 

So, consider $x(0)\in\mathbb{R}^{2}\setminus\mathcal{X}$ , then we
obtain for $x_{1}x_{2}<0$:
\begin{gather*}
\dot{V}=-(\gamma-\text{sign}(x_{1})d)\varepsilon\sqrt{|x_{1})|}-\frac{\bigl|x_{2}\bigr|^{3}}{|x_{1}|}\\
+\varepsilon\frac{x_{2}^{2}}{\sqrt{|x_{1}|}}\left(\frac{1}{2}-\text{sign}(x_{1})\text{sign}(x_{2})\right)+x_{2}d+|x_{2}|(D+\eta)\\
\leq-(\gamma-D)\varepsilon\sqrt{|x_{1}|}-\frac{\bigl|x_{2}\bigr|^{3}}{|x_{1}|}+\varepsilon\frac{3}{2}\frac{x_{2}^{2}}{\sqrt{|x_{1}|}}+|x_{2}|(2D+\eta)\\
\leq-(\gamma-\frac{1}{2}-D-\frac{2}{3\varepsilon}(2D+\eta)^{1.5})\varepsilon\sqrt{|x_{1}|}-\left(\frac{2}{3}-\varepsilon\right)\frac{\bigl|x_{2}\bigr|^{3}}{|x_{1}|},
\end{gather*}
where on the last step the following inequalities have been utilized
(that are obtained from Young's inequality \citep{young1912}): 
\begin{gather*}
\frac{x_{2}^{2}}{\sqrt{|x_{1}|}}\leq\frac{1}{3}\sqrt{|x_{1}|}+\frac{2}{3}\frac{\bigl|x_{2}\bigr|^{3}}{|x_{1}|},\\
(2D+\eta)|x_{2}|\leq\frac{2}{3}(2D+\eta)^{1.5}\sqrt{|x_{1}|}+\frac{1}{3}\frac{\bigl|x_{2}\bigr|^{3}}{|x_{1}|}.
\end{gather*}
Hence, for $\gamma-\frac{1}{2}-D-\frac{2}{3\varepsilon}(2D+\eta)^{1.5}>0$
and $\frac{2}{3}-\varepsilon>0$ the desired property $\dot{V}<0$
while $V>0$ is obtained while $x_{1}x_{2}<0$. Assuming that $\gamma>\frac{1}{2}+D$,
the parameter conditions to check can be reduced to
\[
\frac{2}{3}\frac{(2D+\eta)^{1.5}}{\gamma-\frac{1}{2}-D}<\varepsilon<\min\biggl\{\frac{2}{3},\sqrt{2(\gamma-D-\eta)}\biggr\}.
\]
The interval for the values of $\varepsilon$ is not empty for a sufficiently
small $\eta>0$ provided that the condition (\ref{eq:gamma_new})
is verified. Indeed, this condition implies that $\gamma-D>0.5$ for
any $D>0$, hence, $\min\{\frac{2}{3},\sqrt{2(\gamma-D-\eta)}\}=\frac{2}{3}$
for $\eta\in(0,\frac{1}{3})$ and $\frac{2}{3}\frac{(2D+\eta)^{1.5}}{\gamma-\frac{1}{2}-D}<\frac{2}{3}$
is equivalent to (\ref{eq:gamma_new}) for an infinitesimal $\eta$.

Next, for $x(0)\in\mathbb{R}^{2}\setminus\mathcal{X}$ with $x_{1}x_{2}>0$
we derive:
\begin{gather*}
\dot{V}=(\gamma-D-\eta)\text{sign}(x_{1})\dot{x}_{1}+x_{2}\dot{x_{2}}\\
=(\gamma-D-\eta)\text{sign}(x_{1})x_{2}+x_{2}(d-\gamma\,\textrm{sign}(x_{1}))\\
=-(D+\eta)\text{sign}(x_{1})x_{2}+x_{2}d\leq-\eta|x_{2}|,
\end{gather*}
which implies a strict decay of $V$ again.

Combining the derived estimates for $\dot{V}$ outside of $\mathcal{X}$,
and taking into account that the set $\mathcal{X}\setminus\{0\}$
is visited at a single time instant, the uniform global asymptotic
stability property of the origin for (\ref{eq:2:1})--(\ref{eq:2:3})
follows the conventional results \cite{Lin1996,khalil2002}. 
\end{pf}

As we can conclude, the condition (\ref{eq:gamma_new}) is more restrictive than (\ref{eq:gamma_old}). However, this conservatism does not characterize the properties of the closed-loop system (\ref{eq:2:1})--(\ref{eq:2:3}), since it comes from the Lyapunov function \eqref{eq:V_new} used for the analysis. Moreover, the method used to establish global finite-time convergence to the origin in the previous subsection utilizes (\ref{eq:gamma_old}) only. In addition, it is intuitively clear that even the restriction (\ref{eq:gamma_old}) is related with the analysis of the derivative of \eqref{eq:V_old}, and the finite-time stability can be observed in (\ref{eq:2:1})--(\ref{eq:2:3}) for any $\gamma > D$.

\subsection{Numerical illustration}
\label{sec:4:sub:3}

In the following, we present several numerical illustrations for a perturbed control system \eqref{eq:2:1}--\eqref{eq:2:3}. The numerical simulation framework is the same as the one specified and used in section \ref{sec:3:sub:3}. For disturbances, the upper bound $D=100$ is next assumed.

\begin{figure}[!h] \centering
\includegraphics[width=0.98\columnwidth]{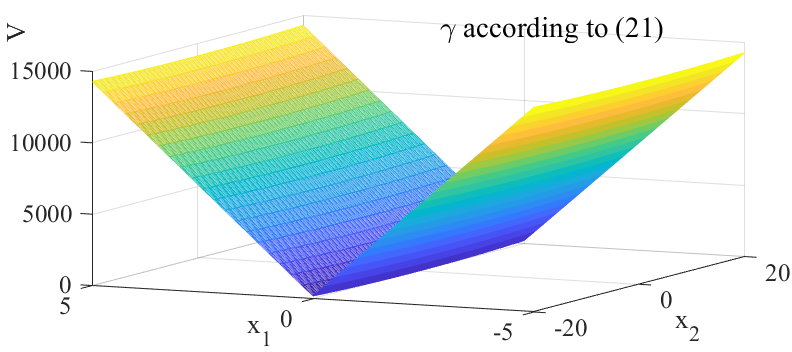}\\[3mm]
\includegraphics[width=0.98\columnwidth]{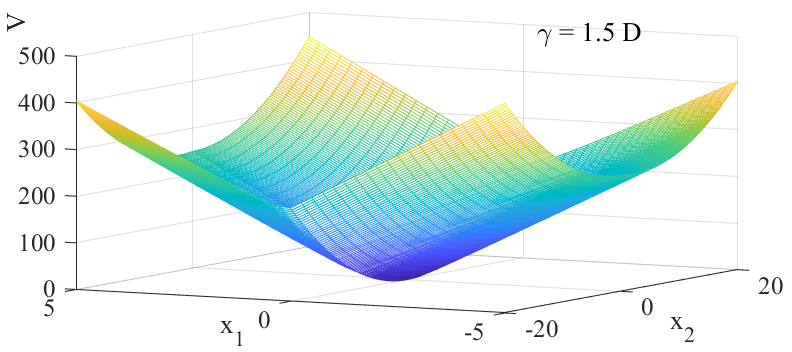}
\caption{Comparison of Lyapunov functions \eqref{eq:V_new} for a conservative $\gamma$ assignment according to \eqref{eq:gamma_new} above, and for intuitively ($\gamma > D$) assigned $\gamma = 1.5 \, D$ below.}
\label{fig:4:1}
\end{figure}
First, we illustrate the shape of the Lyapunov function \eqref{eq:V_new}, once for a conservative gain assignment according to \eqref{eq:gamma_new}, and once for an intuitive $\gamma > D$ assignment  as $\gamma = 1.5 \, D$. Both $V(x)$ maps are depicted in Fig. \ref{fig:4:1}.

\begin{figure}[!h] \centering
\includegraphics[width=0.98\columnwidth]{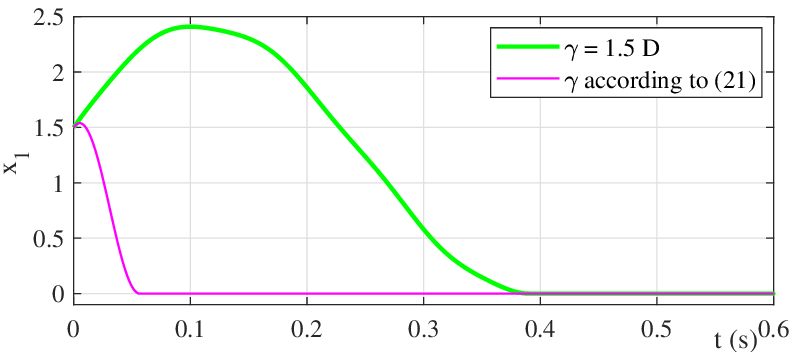}\\[2mm]
\includegraphics[width=0.98\columnwidth]{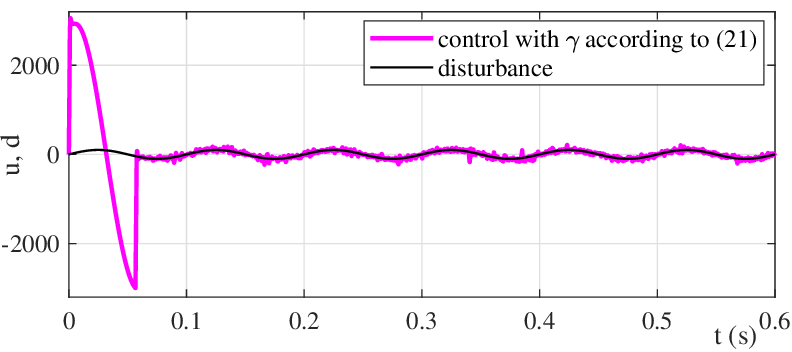}\\[2mm]
\includegraphics[width=0.98\columnwidth]{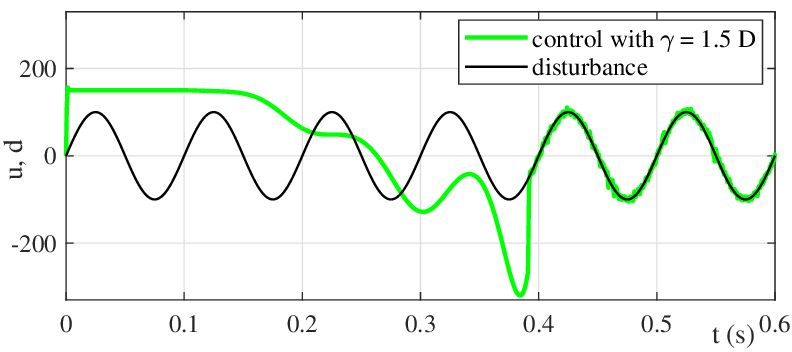}
\caption{Comparison of conservative $\gamma$ assignment according to \eqref{eq:gamma_new} with intuitively ($\gamma > D$) assigned $\gamma = 1.5 \, D$: perturbed output state above and (low-pass filtered) control versus disturbance in the middle and below.}
\label{fig:4:2}
\end{figure}

Next, we compare the output state trajectories and the control values versus disturbance for both assignments of the $\gamma$-gain. The intuitively assigned gain is $\gamma = 1.5 \, D = 150$, while the gain assigned according to \eqref{eq:gamma_new} results in $\gamma = 2929$. The initial state is $x(0) = [1.5, \, 15]$, and the dynamic disturbance is simulated as $d(t) = 100 \sin (20 \pi \,t)$. The disturbance frequency of 10 Hz is assigned by taking into account the convergence time to origin for the set initial values. The perturbed output state trajectories are shown opposite each other in Fig. \ref{fig:4:2} above, while the control versus disturbance values are shown separately for both gain assignments in the middle and below. Note that since the control system proceeds in sliding-mode once being in the origin, the control values are low-pass filtered with 500 Hz cutoff frequency, for the sake of visualization.   

%%%%%%%%%%%%%%%%%%%%%%%%%%%%%%%%%%%%%%%%%%%%%%%%%%%%%%%%%%%%%%%%%%%%%%%%%%%%%%%%
\section{Conclusion}
\label{sec:5}

This paper introduced and analyzed a significant modification of the nonlinear control proposed in \cite{ruderman2024robust}. The modified control law \eqref{eq:2:2}, \eqref{eq:2:3} is well
defined for all states in $\mathbb{R}^2$ and, moreover, the solutions of the closed-loop nonlinear system without perturbations are obtained in the analytic form for a dominant part of the state space. Similar as the initial control (\cite{ruderman2024robust}), the shown control method is finite-time converging without overshooting despite any matched and properly bounded disturbances. These properties come in favor of multiple practical application scenarios and specifications for second-order systems. The uniform global asymptotic stability of the closed-loop system was proven by using the Lyapunov function method. Alongside the developed comprehensive analysis, we provide several numerical illustrations which disclose the dynamic behavior of the control system and discuss the formal conservatism and intuitive criteria for assignment of the single control parameter.

%%%%%%%%%%%%%%%%%%%%%%%%%%%%%%%%%%%%%%%%%%%%%%%%%%%%%%%%%%%%%%%%%%%%%%%%%%%%%%%%
\section*{Acknowledgement}

The first author acknowledges the financial support by NEST (Network for
Energy Sustainable Transition) foundation during the annual sabbatical 
at Polytechnic University of Bari.

\bibliography{references}             % bib file to produce the bibliography
                                                     % with bibtex (preferred)

\end{document}